\documentclass[english]{extarticle}
\usepackage[T1]{fontenc}
\usepackage[latin9]{inputenc}
\usepackage{geometry}
\usepackage{units}
\usepackage{amsmath}
\usepackage{amssymb}
\usepackage{mathdots}
\usepackage{graphicx}
\usepackage[all]{xy}

\makeatletter
\@ifundefined{date}{}{\date{}}

\usepackage[all]{xy}

\usepackage{placeins}

\usepackage{cite}

\AtBeginDocument{
  
}

\makeatother

\usepackage{babel}
\begin{document}
\title{Contents, Contexts, and Basics of Contextuality}
\author{Ehtibar N. Dzhafarov\\$\,$\\Purdue University}
\maketitle
\begin{abstract}
This is a non-technical introduction into theory of contextuality.
More precisely, it presents the basics of a theory of contextuality
called Contextuality-by-Default (CbD). One of the main tenets of CbD
is that the identity of a random variable is determined not only by
its content (that which is measured or responded to) but also by contexts,
systematically recorded conditions under which the variable is observed;
and the variables in different contexts possess no joint distributions.
I explain why this principle has no paradoxical consequences, and
why it does not support the holistic ``everything depends on everything
else'' view. Contextuality is defined as the difference between two
differences: (1) the difference between content-sharing random variables
when taken in isolation, and (2) the difference between the same random
variables when taken within their contexts. Contextuality thus defined
is a special form of context-dependence rather than a synonym for
the latter. The theory applies to any empirical situation describable
in terms of random variables. Deterministic situations are trivially
noncontextual in CbD, but some of them can be described by systems
of epistemic random variables, in which random variability is replaced
with epistemic uncertainty. Mathematically, such systems are treated
as if they were ordinary systems of random variables. 
\end{abstract}

\section{Contents, contexts, and random variables}

The word \emph{contextuality} is used widely, usually as a synonym
of \emph{context-dependence}. Here, however, contextuality is taken
to mean a special form of context-dependence, as explained below.
Historically, this notion is derived from two independent lines of
research: in quantum physics, from studies of existence or nonexistence
of the so-called hidden variable models with context-independent mapping
\cite{Bell1964,CHSH1969,CH1974,Fine1982JMP,KurzynskiCabelloKaszlikowski(2014),KochenSpecker(1967),Cabello(2008),Cabelloetal(1996),KCBS2008,Bell1966},\footnote{Here, I mix together the early studies of nonlocality and those of
contextuality in the narrow sense, related to the Kochen-Specker theorem
\cite{KochenSpecker(1967)}. Both are special cases of contextuality.} and in psychology, from studies of the so-called selective influences
\cite{Dzhafarov(2003),DzhafarovGluhovsky(2006),DzhafarovKujala(2010),DzhafarovKujala(2016)NHMP,Sternberg1969,Townsend1984,KujalaDzhafarov(2008),ZhangDzhafarov(2015)}.
The two lines of research merged relatively recently, in the 2010's
\cite{DzhafarovKujala(2012)Selectivity,DzhafarovKujala(2012)Quantum,DzhafarovKujala(2013)AllPossible,DzhafarovKujala(2014)AQualifiesKolmogorov,DzhafarovKujala(2013)OrderDistance,DzhafarovKujala(2014)TopicsCogSci},
to form an abstract mathematical theory, Contextuality-by-Default
(CbD), with multidisciplinary applications \cite{Dzh2017Nothing,Dzhafarov(2016),DzhafarovKujala(2016)Context-Content,DzhafarovKujala(2017)LNCS-2,DzhafarovKujalaCervantes(2016)LNCS,DzhafarovKujalaLarsson(2015),DzhafarovZhangKujala(2015)Isthere,DzhCerKuj2017,DzhKuj2017Fortsch-1,DzhKuj2018,KujalaDzhafarov(2016)Proof,ZhangDzhafarov(2016),Basievaetall2018,bacciagaluppi2,CD2020Impossible,CervantesDzhafarov2017,CervDzh2018,CervDzhaf2019,deBarrosDzhafarovetal(2016),DKC2020,DKC2020Erratum,DKC2020DEpistemic,Dzhafarov2019,Dzhafarov2020,DzhKuj2020,DKconversations,DzhafarovKujalaCervantesZhangJones(2016),KujalaDzhafarovLar(2015),KujDzh2019,Jones2019,DK2014Scripta-1,CervantesDzhafarov(2016),KujalaDzhafarov2015}.\footnote{The theory has been revised in two ways since 2016, the changes being
presented in Refs. \cite{DzhafarovKujala(2017)LNCS-2,DzhCerKuj2017}.}

The example I will use to introduce the notion of contextuality reflects
the fact that even as I write these lines the world is being ravaged
by the Covid-19 pandemic, forcing lockdowns and curtailing travel. 

Suppose we ask a randomly chosen person two questions: 

\[
\begin{array}{|l|}
\hline q_{1}:\textnormal{would you like to take an overseas vacation this summer}?\\
\hline q_{2}:\textnormal{are you wary of contracting Covid-19}?
\\\hline \end{array}
\]
Suppose also we ask these questions in two orders:
\[
\begin{array}{|l|}
\hline c^{1}:\textnormal{first }q_{1}\textnormal{ then }q_{2}\\
\hline c^{2}:\textnormal{first }q_{2}\textnormal{ then }q_{1}
\\\hline \end{array}
\]
To each of the two questions, the person can respond in one of two
ways: Yes or No. And since we are choosing people to ask our questions
randomly, we cannot determine the answer in advance. We assume therefore
that the answers can be represented by \emph{random variables}. A
random variable is characterized by its \emph{identity} (as explained
shortly) and its \emph{distribution}: in this case, the distribution
means responses Yes and No together with their probabilities of occurrence.\footnote{I set aside the intriguing issue of whether responses Yes and No may
be indeterministic but not assignable probabilities.}

One can summarize this imaginary experiment in the form of the following
\emph{system of random variables}:

\begin{equation}
\begin{array}{|c|c||c|}
\hline R_{1}^{1} & R_{2}^{1} & c^{1}=q_{1}\rightarrow q_{2}\\
\hline R_{1}^{2} & R_{2}^{2} & c^{2}=q_{2}\rightarrow q_{1}\\
\hline\hline q_{1}=\textnormal{"}\textnormal{vacation?"} & q_{2}=\textnormal{"}\textnormal{Covid-19?"} & \textnormal{system }\mathcal{C}_{2(a)}
\\\hline \end{array}\:.\label{eq:matrix1}
\end{equation}
This is the simplest system that can exhibit contextuality (as defined
below). The random variables representing responses to questions are
denoted by $R$ with subscripts and superscripts determining its identity.
The subscript of a random variable in the system refers to the question
this random variable answers: e.g., $R_{1}^{1}$ and $R_{1}^{2}$
both answer the question $q_{1}$. The superscript refers to the \emph{context}
of the random variable, the circumstances under which it is recorded.
In the example the context is the order in which the two questions
are being asked. Thus, $R_{2}^{1}$ answers question $q_{2}$ when
this question is asked second, whereas $R_{2}^{2}$ answers the same
question when it is is asked first.

The question a random variable answers is generically referred to
as this variable's \emph{content}. Contents can always be thought
of as having the logical function of questions, but in many cases
other than in our example they are not questions in the colloquial
meaning. Thus, a $q$ may be one's choice of a physical object to
measure, say, a stone to weigh, in which case the stone will be the
content of the random variable $R_{q}^{c}$ representing the outcome
of weighing it (in some context $c$). Of course, logically, this
$R_{q}^{c}$ answers the question of how heavy the stone is, and $q$
can be taken to stand for this question. 

Returning to our example, each variable $R_{q}^{c}$ in our set of
four variables is identified by its content ($q=q_{1}$ or $q=q_{2}$)
and by its context ($c=c^{1}$ or $c=c^{2}$). It is this double-identification
that imposes a structure on this set, rendering it a \emph{system}
(specifically, a \emph{content-context system}) of random variables.
There may be other variable circumstances under which our questions
are asked, such as when and where the questions were asked, in what
tone of voice, or how high the solar activity was when they were asked.
However, it is a legitimate choice not to take such concomitant circumstances
into account, to ignore them. If we do not, which is a legitimate
choice too, our contexts will have to be redefined, yielding a different
system, with more than just four random variables. The legitimacy
of ignoring all but a select set of contexts is an important aspect
of contextuality analysis, as we will see later. 

\begin{figure}
\[
\vcenter{\xymatrix@C=1cm{ &  & {\scriptstyle R_{1}^{1}}\ar@{-}[r]_{\textnormal{context}c^{1}} & {\scriptstyle R_{2}^{1}}\ar@{.}[dr]_{\textnormal{content }q_{2}}\\
 & {\scriptstyle R_{1}^{n}}\ar@{.}[ur]_{\textnormal{content}q_{1}} &  &  & {\scriptstyle R_{2}^{2}}\ar@{-}[dr]_{\textnormal{context}c^{2}}\\
{\scriptstyle R_{n}^{n}}\ar@{-}[ur]_{\textnormal{context}c^{n}} &  &  &  &  & {\scriptstyle R_{3}^{2}}\ar@{.}[dl]_{\textnormal{content}q_{3}}\\
 & \ddots\ar@{.}[ul]_{\textnormal{content}q_{n}} &  &  & \iddots\ar@{.}[dl]_{\textnormal{content}q_{i}}\\
 &  & {\scriptstyle R_{i+1}^{i}}\ar@{.}[lu]_{\textnormal{content }q_{i+1}} & {\scriptstyle R_{i}^{i}}\ar@{-}[l]_{\textnormal{context }c^{i}}
}
}
\]

\caption{\label{fig:A-cyclic-system}A cyclic system of rank $n$.}

\end{figure}
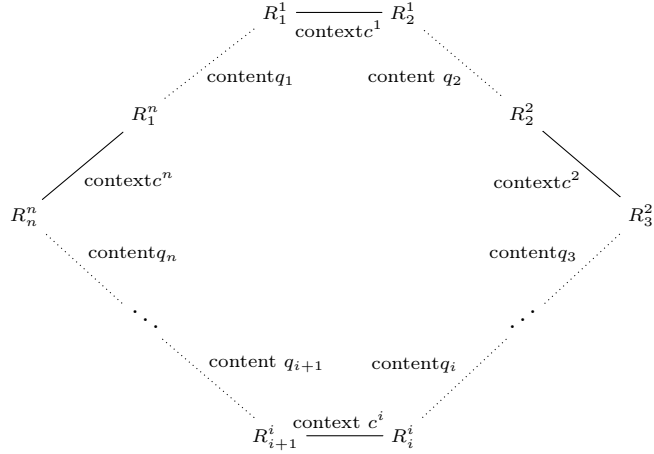

The reason I denote our system $\mathcal{C}_{2(a)}$ is that it is
a specific example (the specificity being indicated by index $a$)
of a \emph{cyclic system of rank }2, denoted $\mathcal{C}_{2}$. More
generally, cyclic systems of rank $n$, denoted $\mathcal{C}_{n}$,
are characterized by the arrangement of $n$ contents, $n$ contexts,
and $2n$ random variables shown in Figure \ref{fig:A-cyclic-system}.

A system of the $\mathcal{C}_{2}$-type is the smallest such system
(not counting the degenerate system consisting of $R_{1}^{1}$ alone):
\[
\vcenter{\xymatrix@C=1cm{{\scriptstyle R_{1}^{1}}\ar@{-}[rrr]_{\textnormal{context}c^{1}} &  &  & {\scriptstyle R_{2}^{1}}\ar@{.}[dd]_{\textnormal{content }q_{2}}\\
\\
{\scriptstyle R_{1}^{2}}\ar@{.}[uu]_{\textnormal{content }q_{1}} &  &  & {\scriptstyle R_{2}^{2}}\ar@{-}[lll]_{\textnormal{context }c^{2}}
}
}\:.
\]

What else do we know of our random variables? First of all, the two
variables within a context, $\left(R_{1}^{1},R_{2}^{1}\right)$, or
$\left(R_{1}^{2},R_{2}^{2}\right)$, are \emph{jointly distributed}.
By the virtue of being responses of one and same person, the values
of these random variables come in pairs. So it is meaningful to ask
what the probabilities are for each of the joint events
\[
\begin{array}{ccc}
R_{1}^{1}=+1 & \textnormal{and} & R_{2}^{1}=+1,\\
R_{1}^{1}=+1 & \textnormal{and} & R_{2}^{1}=-1,\\
R_{1}^{1}=-1 & \textnormal{and} & R_{2}^{1}=+1,\\
R_{1}^{1}=-1 & \textnormal{and} & R_{2}^{1}=-1,
\end{array}
\]
where $+1$ and $-1$ encode the answers Yes and No, respectively.
One can meaningfully speak of correlations between the variables in
the same context, probability that they have the same value, etc.

By contrast, different contexts, in our case the two orders in which
the questions are asked, are mutually exclusive. When asked two questions,
a given person can only be asked them in one order. Respondents represented
by $R_{1}^{1}$ answer question $q_{1}$ asked first, before $q_{2}$,
whereas the respondents represented by $R_{1}^{2}$ answer question
$q_{1}$ asked second, after $q_{2}$. Clearly, these are different
sets of respondents, and one would not know how to pair them. It is
meaningless to ask, e.g., what the probability of 
\[
\begin{array}{ccc}
R_{1}^{1}=+1 & \textnormal{and} & R_{1}^{2}=+1\end{array}
\]
may be. Random variables in different contexts are \emph{stochastically
unrelated}. 

\section{\label{sec:Intuition-of-(non)contextuality}Intuition of (non)contextuality}

Having established these basic facts, let us consider now the two
random variables with content $q_{1}$, and let us make at first the
(unrealistic) assumption that their distributions are the same in
both contexts, $c^{1}$ and $c^{2}$: 

\begin{equation}
\begin{array}{c|c}
value & probability\\
\hline R_{1}^{1}=+1 & a\\
\hline R_{1}^{1}=-1 & 1-a
\end{array}\textnormal{ and }\begin{array}{c|c}
value & probability\\
\hline R_{1}^{2}=+1 & a\\
\hline R_{1}^{2}=-1 & 1-a
\end{array}\:.\label{eq:consist connect 1}
\end{equation}
If we consider the variables $R_{1}^{1}$ and $R_{1}^{2}$ in isolation
from their contexts (i.e., disregarding the other two random variables),
then we can view them as simply one and the same random variable.
In other words, the subsystem

\[
\begin{array}{|c||c|}
\hline R_{1}^{1} & c^{1}=q_{1}\rightarrow q_{2}\\
\hline R_{1}^{2} & c^{2}=q_{2}\rightarrow q_{1}\\
\hline\hline q_{1}=\textnormal{"}\textnormal{vacation?}" & \mathcal{C}_{2(a)}/\textnormal{only }q_{1}
\\\hline \end{array}
\]
appears to be replaceable with just
\[
\begin{array}{|c|}
\hline R_{1}\\
\hline\hline q_{1}=\textnormal{"}\textnormal{vacation?}"
\\\hline \end{array}\:,
\]
with contexts being superfluous. 

Analogously, if the distributions of the two random variables with
content $q_{2}$ are assumed to be the same,
\begin{equation}
\begin{array}{c|c|c}
value & R_{2}^{1}=+1 & R_{2}^{1}=-1\\
\hline probability & b & 1-b
\end{array}\textnormal{ and }\begin{array}{c|c|c}
value & R_{2}^{2}=+1 & R_{2}^{2}=-1\\
\hline probability & b & 1-b
\end{array}\:,\label{eq:consist connect 2}
\end{equation}
and if we consider them in isolation from their contexts, the subsystem

\[
\begin{array}{|c||c|}
\hline R_{2}^{1} & c^{1}=q_{1}\rightarrow q_{2}\\
\hline R_{2}^{2} & c^{2}=q_{2}\rightarrow q_{1}\\
\hline\hline q_{2}=\textnormal{"}\textnormal{Covid-19?}" & \mathcal{C}_{2(a)}/\textnormal{only }q_{2}
\\\hline \end{array}
\]
appears to be replaceable with

\[
\begin{array}{|c|}
\hline R_{2}\\
\hline\hline q_{2}=\textnormal{"}\textnormal{Covid-19?}"
\\\hline \end{array}\:.
\]
It is tempting now to say: we have only two random variables, $R_{1}$
and $R_{2}$, whatever their contexts. But a given pair of random
variables can only have one \emph{joint distribution}, this distribution
cannot be somehow different in different contexts. We should predict
therefore, that if the probabilities in system $\mathcal{C}_{2(a)}$
are
\[
\Pr\left[R_{1}^{1}=+1,R_{2}^{1}=+1\right]=r_{1}\textnormal{ and }\Pr\left[R_{1}^{2}=+1,R_{2}^{2}=+1\right]=r_{2},
\]
then 
\[
r_{1}=r_{2}.
\]

Suppose, however, that this is shown to be empirically false, that
in fact $r_{1}>r_{2}$. For instance, assuming $0<a<b$, suppose that
the joint distributions in the two contexts of system $\mathcal{C}_{2(a)}$
are 
\begin{equation}
\begin{array}{c|c|c|c}
\textnormal{context }c^{1} & R_{2}^{1}=+1 & R_{2}^{1}=-1\\
\hline R_{1}^{1}=+1 & r_{1}=a & 0 & a\\
\hline R_{1}^{1}=-1 & b-a & 1-b & 1-a\\
\hline  & b & 1-b
\end{array}\label{eq:context1 of system1}
\end{equation}
and

\begin{equation}
\begin{array}{c|c|c|c}
\textnormal{context }c^{2} & R_{2}^{2}=+1 & R_{2}^{2}=-1\\
\hline R_{1}^{2}=+1 & r_{2}=0 & a & a\\
\hline R_{1}^{2}=-1 & b & 1-a-b & 1-a\\
\hline  & b & 1-b
\end{array}\:.\label{eq:context2 of system1}
\end{equation}
Clearly, we have then a \emph{reductio ad absurdum} proof that the
assumption we have made is wrong, the assumption being that we can
drop contexts in $R_{1}^{1}$ and $R_{1}^{2}$ (as well as in $R_{2}^{1}$
and $R_{2}^{2}$), and that we can therefore treat them as one and
the same random variable $R_{1}$ (respectively, $R_{2}$). This is
the simplest case when we can say that a system of random variables,
here, the system $\mathcal{C}_{2(a)}$, is \emph{contextual}. 

This understanding of contextuality can be extended to more complex
systems. However, it is far from being general enough. It only applies
to \emph{consistently connected} systems, those in which any two variables
with the same content are identically distributed.\footnote{The term ``consistent connectedness'' is due to the fact that in
CbD the content-sharing random variables are said to form \emph{connections}
(between contexts). In quantum physics consistent connectedness is
referred to by such terms as lack of signaling, lack of disturbance,
parameter invariance, etc.} This assumption is often unrealistic. Specifically, it is a well-established
empirical fact that the individual distributions of the responses
to two questions do depend on their order \cite{Moore}. Besides,
this is highly intuitive in our example. If one is asked about an
overseas vacation first, the probability of saying ``Yes, I would
like to take an overseas vacation'' may be higher than when this
question is asked second, after the respondent has been reminded about
the dangers of the pandemic. 

In order to generalize the notion of contextuality to arbitrary systems,
we need to develop answers to the following two questions:
\begin{description}
\item [{A:}] For any two random variables sharing a content, how different
are they when taken in isolation from their contexts? 
\item [{B:}] Can these differences be preserved when all pairs of content-sharing
variables are taken within their contexts (i.e., taking into account
their joint distributions with other random variables in their contexts)?
\end{description}
For our system $\mathcal{C}_{2(a)}$ with the within-context joint
distributions given by (\ref{eq:context1 of system1}) and (\ref{eq:context2 of system1}),
our informal answer to question \textbf{A} was that two random variables
with the same content (i.e., $R_{1}^{1}$ and $R_{1}^{2}$ or $R_{2}^{1}$
and $R_{2}^{2}$) are not different at all when taken in isolation.
The informal answer to question \textbf{B}, however, was that in these
two pairs (or at least in one of them) the random variables are not
the same when taken in relation to other random variables in their
respective contexts. One can say therefore that 
\begin{quote}
\emph{the contexts make $R_{1}^{1}$ and $R_{1}^{2}$ (and/or $R_{2}^{1}$
and $R_{2}^{2}$) more dissimilar than when they are taken without
their contexts}. 
\end{quote}
This is the intuition we will use to construct a general definition
of contextuality.

\section{Making it rigorous: Couplings}

First, we have to agree on how to measure the difference between two
random variables that are not jointly distributed, like $R_{1}^{1}$
and $R_{1}^{2}$. Denote these random variables $X$ and $Y$, both
dichotomous ($\pm1$), with
\[
\Pr\left[X=+1\right]=u\textnormal{ and }\Pr\left[Y=+1\right]=v.
\]
 Consider all possible pairs of \emph{jointly distributed }variables
$\left(X',Y'\right)$ such that
\[
X'\overset{dist}{=}X,Y'\overset{dist}{=}Y,
\]
where $\overset{dist}{=}$ stands for ``has the same distribution
as.'' Any such pair $\left(X',Y'\right)$ is called a \emph{coupling}
of $X$ and $Y$. For obvious reasons, two couplings of $X$ and $Y$
having the same joint distribution are not distinguished. 

Now, for each coupling $\left(X',Y'\right)$ one can compute the probability
with which $X'\not=Y'$ (recall that the probability of $X\not=Y$
is undefined, we do need couplings to make this inequality a meaningful
event). It is easy to see that among the couplings $\left(X',Y'\right)$
there is one and only one for which this probability is minimal. This
coupling is defined by the joint distribution
\begin{equation}
\begin{array}{c|c|c|c}
 & Y'=+1 & Y'=-1\\
\hline X'=+1 & \min\left(u,v\right) & u-\min\left(u,v\right) & u\\
\hline X'=-1 & v-\min\left(u,v\right) & \min\left(1-u,1-v\right) & 1-u\\
\hline  & v & 1-v
\end{array}\:,\label{eq:general max coupling}
\end{equation}
and the minimal probability in question is obtained as
\[
\left(u-\min\left(u,v\right)\right)+\left(v-\min\left(u,v\right)\right)=\left|u-v\right|.
\]
This probability is a natural measure of difference between the random
variables $X$ and $Y$:\footnote{It is a special case of the so-called \emph{total variation distance},
except that it is usually defined between two probability distributions,
while I use it here as a measure of difference (formally, a \emph{pseudometric})
between stochastically unrelated random variables.}
\begin{equation}
\delta\left(X,Y\right)=\min_{\begin{array}{c}
\textnormal{all couplings}\\
\left(X',Y'\right)\textnormal{ of }\ensuremath{X\textnormal{ and }Y}
\end{array}\textnormal{ }}\Pr\left[X'\not=Y'\right]=\left|u-v\right|.\label{eq:general min difference}
\end{equation}
If $X$ and $Y$ are identically distributed, i.e. $u=v$, the joint
distribution of $X'$ and $Y'$ can be chosen as 
\[
\begin{array}{c|c|c|c}
\textnormal{context }c^{1} & Y=+1 & Y=-1\\
\hline X=+1 & u & 0 & u\\
\hline X=-1 & 0 & 1-u & 1-u\\
\hline  & u & 1-u
\end{array},
\]
yielding 
\[
\delta\left(X,Y\right)=\min_{\begin{array}{c}
\textnormal{all couplings}\\
\left(X',Y'\right)\textnormal{ of }\ensuremath{X\textnormal{ and }Y}
\end{array}\textnormal{ }}\Pr\left[X'\not=Y'\right]=0.
\]

Let us apply this to our example, in order to formalize the intuition
behind our saying earlier that two identically distributed random
variables, taken in isolation, can be viewed as being ``the same.''
For $R_{1}^{1}$ and $R_{1}^{2}$ in (\ref{eq:consist connect 1}),
\[
\delta\left(R_{1}^{1},R_{1}^{2}\right)=\min_{\begin{array}{c}
\textnormal{all couplings}\\
\left(S_{1}^{1},S_{1}^{2}\right)\textnormal{ of }\ensuremath{R_{1}^{1}\textnormal{ and }R_{1}^{2}}
\end{array}\textnormal{ }}\Pr\left[S_{1}^{1}\not=S_{1}^{2}\right]=0,
\]
and, analogously, for $R_{2}^{1}$ and $R_{2}^{2}$ in (\ref{eq:consist connect 2}),
\[
\delta\left(R_{2}^{1},R_{2}^{2}\right)=\min_{\begin{array}{c}
\textnormal{all couplings}\\
\left(S_{2}^{1},S_{2}^{2}\right)\textnormal{ of }\ensuremath{R_{2}^{1}\textnormal{ and }R_{2}^{2}}
\end{array}\textnormal{ }}\Pr\left[S_{2}^{1}\not=S_{2}^{2}\right]=0.
\]

\section{Making it rigorous: Contextuality}

What is then the rigorous way of establishing that these differences
cannot both be zero when considered within their contexts? For this,
we need to extend the notion of a coupling to an entire system. A
coupling of our system $\mathcal{C}_{2(a)}$ is a set of \emph{corresponding
jointly distributed} random variables 
\begin{equation}
\begin{array}{|c|c|}
\hline S_{1}^{1} & S_{2}^{1}\\
\hline S_{1}^{2} & S_{2}^{2}
\\\hline \end{array}\label{eq:coupling1}
\end{equation}
such that 
\begin{equation}
\left(S_{1}^{1},S_{2}^{1}\right)\overset{dist}{=}\left(R_{1}^{1},R_{2}^{1}\right),\left(S_{1}^{2},S_{2}^{2}\right)\overset{dist}{=}\left(R_{1}^{2},R_{2}^{2}\right).\label{eq:joints1}
\end{equation}
In other words, the distributions within contexts, (\ref{eq:context1 of system1})
and (\ref{eq:context2 of system1}), remain intact when we replace
the $R$-variables with the corresponding $S$-variables, 
\begin{equation}
\begin{array}{c|c|c|c}
 & S_{2}^{1}=+1 & S_{2}^{1}=-1\\
\hline S_{1}^{1}=+1 & a & 0 & a\\
\hline S_{1}^{1}=-1 & b-a & 1-b & 1-a\\
\hline  & b & 1-b
\end{array}\textnormal{ and }\begin{array}{c|c|c|c}
 & S_{2}^{2}=+1 & S_{2}^{2}=-1\\
\hline S_{1}^{2}=+1 & 0 & a & a\\
\hline S_{1}^{2}=-1 & b & 1-a-b & 1-a\\
\hline  & b & 1-b
\end{array}\:.\label{eq:same probs}
\end{equation}

Such couplings always exist, not only for our example, but for any
other system of random variables. Generally, there is an infinity
of couplings for a given system.\footnote{One need not have separate definitions of couplings for pairs of random
variables and for systems. In general, given any set of random variables
$\mathfrak{R}$, its coupling is a set of random variables $S$, in
a one-to-one correspondence with $\mathfrak{R}$, such that the corresponding
variables in $\mathfrak{R}$ and $S$ have the same distribution,
and all variables in $S$ are jointly distributed. To apply this definition
to $\mathfrak{R}$ representing a system of random variables one considers
all variables within a given context as a single element of $\mathfrak{R}$.
In our example, (\ref{eq:coupling1}) is a coupling of two stochastically
unrelated random variables, $\left(R_{1}^{1},R_{2}^{1}\right)$ and
$\left(R_{1}^{2},R_{2}^{2}\right)$.} Thus, to construct a coupling for system $\mathcal{C}_{2(a)}$, one
has to assign probabilities to all quadruples of joint events,
\[
\begin{array}{cccc|c|c}
S_{1}^{1} & S_{2}^{1} & S_{1}^{2} & S_{2}^{2} &  & probability\\
\hline +1 & +1 & +1 & +1 &  & p_{++++}\\
+1 & +1 & +1 & -1 &  & p_{+++-}\\
\vdots & \vdots & \vdots & \vdots &  & \vdots\\
-1 & -1 & -1 & -1 &  & p_{----}
\end{array}
\]
so that the appropriately chosen subsets of these probabilities sum
to the joint probabilities shown in (\ref{eq:same probs}):
\[
\begin{array}{c}
p_{++++}+p_{+++-}+p_{++-+}+p_{++--}=\Pr\left[S_{1}^{1}=+1,S_{2}^{1}=+1\right]=a,\\
p_{+-++}+p_{+-+-}+p_{+--+}+p_{+---}=\Pr\left[S_{1}^{1}=+1,S_{2}^{1}=-1\right]=0,\\
p_{++++}+p_{+-++}+p_{-+++}+p_{--++}=\Pr\left[S_{1}^{2}=+1,S_{2}^{2}=+1\right]=0,\\
etc.
\end{array}
\]
This is a system of seven independent linear equations with 16 unknown
$p$-probabilities, subject to the additional constraint that all
probabilities must be nonnegative. It can be shown that this linear
programming problem always has solutions, and infinitely many of them
at that, unless one of the probabilities $a$ and $b$ equals 1 or
0 (in which case the solution is unique).

Unlike in system $\mathcal{C}_{2(a)}$ itself, in any coupling (\ref{eq:coupling1})
of this system the random variables have joint distributions across
the contexts. In particular, $\left(S_{1}^{1},S_{1}^{2}\right)$ is
a jointly distributed pair. Since from (\ref{eq:joints1}) we know
that
\[
S_{1}^{1}\overset{dist}{=}R_{1}^{1}\textnormal{ and }S_{1}^{2}\overset{dist}{=}R_{1}^{2},
\]
$\left(S_{1}^{1},S_{1}^{2}\right)$ is a coupling of $R_{1}^{1}$
and $R_{2}^{1}$. Similarly, $\left(S_{2}^{1},S_{2}^{2}\right)$ is
a coupling of $R_{2}^{1}$ and $R_{2}^{2}$. We ask now: what are
the possible values of
\[
\Pr\left[S_{1}^{1}\not=S_{1}^{2}\right]\textnormal{ and }\Pr\left[S_{2}^{1}\not=S_{2}^{2}\right]
\]
across all possible couplings (\ref{eq:coupling1}) of the entire
system $\mathcal{C}_{2(a)}$? Consider two cases.

\textbf{Case 1.} In some of the couplings (\ref{eq:coupling1}), 
\[
\Pr\left[S_{1}^{1}\not=S_{1}^{2}\right]=0\textnormal{ and }\Pr\left[S_{2}^{1}\not=S_{2}^{2}\right]=0.
\]
We can say then that both $\delta\left(R_{1}^{1},R_{1}^{2}\right)$
and $\delta\left(R_{1}^{1},R_{1}^{2}\right)$ preserve their individual
(in-isolation) values when considered within the system. The system
$\mathcal{C}_{2(a)}$ is then considered \emph{noncontextual}. 

\textbf{Case 2.} In all couplings (\ref{eq:coupling1}), at least
one of the values 
\[
\Pr\left[S_{1}^{1}\not=S_{1}^{2}\right]\textnormal{ and }\Pr\left[S_{2}^{1}\not=S_{2}^{2}\right]
\]
is greater than zero. That is, when considered within the system,
$\delta\left(R_{1}^{1},R_{1}^{2}\right)$ and $\delta\left(R_{1}^{1},R_{1}^{2}\right)$
cannot both be zero. Intuitively, the contexts ``force'' either
$R_{1}^{1}$ and $R_{1}^{2}$ or $R_{2}^{1}$ and $R_{2}^{2}$ (or
both) to be more dissimilar than when taken in isolation. The system
$\mathcal{C}_{2(a)}$ is then considered \emph{contextual}. 

We can quantify the degree of contextuality in the system in the following
way. We know that
\begin{multline*}
\delta\left(R_{1}^{1},R_{1}^{2}\right)+\delta\left(R_{2}^{1},R_{2}^{2}\right)\\
\\
=\min_{\begin{array}{c}
\textnormal{all couplings}\\
\left(S_{1}^{1},S_{1}^{2}\right)\textnormal{ of }\ensuremath{R_{1}^{1}\textnormal{ and }R_{1}^{2}}
\end{array}\textnormal{ }}\left(\Pr\left[S_{1}^{1}\not=S_{1}^{2}\right]\right)+\min_{\begin{array}{c}
\textnormal{all couplings}\\
\left(S_{2}^{1},S_{2}^{2}\right)\textnormal{ of }\ensuremath{R_{2}^{1}\textnormal{ and }R_{2}^{2}}
\end{array}\textnormal{ }}\left(\Pr\left[S_{2}^{1}\not=S_{2}^{2}\right]\right)=0.
\end{multline*}
This quantity is compared to
\begin{multline*}
\delta\left(\left(R_{1}^{1},R_{1}^{2}\right),\left(R_{2}^{1},R_{2}^{2}\right)\right)=\min_{\begin{array}{c}
\textnormal{all couplings}\\
\left(S_{1}^{1},S_{2}^{1},S_{1}^{2},S_{2}^{2}\right)\textnormal{ of system }\mathcal{C}_{2(a)}
\end{array}\textnormal{ }}\left(\Pr\left[S_{1}^{1}\not=S_{1}^{2}\right]+\Pr\left[S_{2}^{1}\not=S_{2}^{2}\right]\right),
\end{multline*}
which can be interpreted as the total of the pairwise differences
between same-content variables within the system. The system is contextual
if this quantity is greater than zero, and this quantity can be taken
as a measure of the degree of contextuality. This is by far not the
only possible measure, but it is arguably the simplest one within
the conceptual framework of CbD.

\section{Generalizing to arbitrary systems}

Consider now a realistic version of our example, when
\[
\begin{array}{c}
\Pr\left[R_{1}^{1}=+1\right]=a_{1},\Pr\left[R_{2}^{1}=+1\right]=b_{1},\\
\Pr\left[R_{1}^{2}=+1\right]=a_{2},\Pr\left[R_{2}^{2}=+1\right]=b_{2},
\end{array}
\]
with $a_{1}$ allowed to be different from $a_{2}$, and $b_{1}$
from $b_{2}$. The within-context joint distributions then generally
look like this: 
\begin{equation}
\begin{array}{c|c|c|c}
\textnormal{context }c^{1} & R_{2}^{1}=+1 & R_{2}^{1}=-1\\
\hline R_{1}^{1}=+1 & r_{1} & a_{1}-r_{1} & a_{1}\\
\hline R_{1}^{1}=-1 & b_{1}-r_{1} & 1-a_{1}-b_{1}+r_{1} & 1-a_{1}\\
\hline  & b_{1} & 1-b_{1}
\end{array}\label{eq:context1 of system1-1}
\end{equation}
and

\begin{equation}
\begin{array}{c|c|c|c}
\textnormal{context }c^{2} & R_{2}^{2}=+1 & R_{2}^{2}=-1\\
\hline R_{1}^{2}=+1 & r_{2} & a_{2}-r_{2} & a_{2}\\
\hline R_{1}^{2}=-1 & b_{2}-r_{2} & 1-a_{2}-b_{2}+r_{2} & 1-a_{2}\\
\hline  & b_{2} & 1-b_{2}
\end{array}\:.\label{eq:context2 of system1-1}
\end{equation}

Let us call the system in (\ref{eq:matrix1}) with these within-context
distributions $C_{2(b)}$. We clearly have context-dependence now
(unless the two joint distributions are identical), but can we also
say that the system is contextual? If we follow the logic of the definition
of contextuality as it was presented above, for consistently connected
systems, the answer cannot automatically be affirmative. The logic
in question requires that we answer the questions \textbf{A} and \textbf{B}
formulated in Section \ref{sec:Intuition-of-(non)contextuality}.
By now we have all necessary conceptual tools for this.

To answer \textbf{A} we look at all possible couplings $\left(S_{1}^{1},S_{1}^{2}\right)$
and $\left(S_{2}^{1},S_{2}^{2}\right)$ of the content-sharing pairs
$\left\{ R_{1}^{1},R_{1}^{2}\right\} $ and $\left\{ R_{2}^{1},R_{2}^{2}\right\} $,
respectively, and determine 
\[
\delta\left(R_{1}^{1},R_{1}^{2}\right)=\min_{\begin{array}{c}
\textnormal{all couplings}\\
\left(S_{1}^{1},S_{1}^{2}\right)\textnormal{ of }\left\{ R_{1}^{1},R_{1}^{2}\right\} 
\end{array}\textnormal{ }}\Pr\left[S_{1}^{1}\not=S_{1}^{2}\right],
\]
and
\[
\delta\left(R_{2}^{1},R_{2}^{2}\right)=\min_{\begin{array}{c}
\textnormal{all couplings}\\
\left(S_{2}^{1},S_{2}^{2}\right)\textnormal{ of }\left\{ R_{2}^{1},R_{2}^{2}\right\} 
\end{array}\textnormal{ }}\Pr\left[S_{2}^{1}\not=S_{2}^{2}\right].
\]
To answer \textbf{B}, we look at all possible couplings
\[
\begin{array}{|c|c|}
\hline S_{1}^{1} & S_{2}^{1}\\
\hline S_{1}^{2} & S_{2}^{2}
\\\hline \end{array}
\]
of the entire system $\mathcal{C}_{2(b)}$, and determine if we can
find couplings in which
\[
\Pr\left[S_{1}^{1}\not=S_{1}^{2}\right]=\delta\left(R_{1}^{1},R_{1}^{2}\right)
\]
and
\[
\Pr\left[S_{2}^{1}\not=S_{2}^{2}\right]=\delta\left(R_{2}^{1},R_{2}^{2}\right).
\]
If such couplings exist, we say that the system is noncontextual,
even if it exhibits context-dependence in the form of \emph{inconsistent
connectedness}. 

Recall that consistently connected systems are those in which any
two variables with the same content are identically distributed, as
it was in our initial (unrealistic) example. For such systems $\delta\left(R_{1}^{1},R_{1}^{2}\right)=0$
and $\delta\left(R_{2}^{1},R_{2}^{2}\right)=0$. However, if 
\[
R_{1}^{1}\overset{dist}{\not=}R_{1}^{2},
\]
then $\delta\left(R_{1}^{1},R_{1}^{2}\right)>0$, and analogously
for $\delta\left(R_{2}^{1},R_{2}^{2}\right)$. In fact, we know from
(\ref{eq:general max coupling}) and (\ref{eq:general min difference})
that if the within-context distributions in the system are as in (\ref{eq:context1 of system1-1})
and (\ref{eq:context2 of system1-1}), then
\[
\delta\left(R_{1}^{1},R_{1}^{2}\right)=\left|a_{1}-a_{2}\right|,\delta\left(R_{2}^{1},R_{2}^{2}\right)=\left|b_{1}-b_{2}\right|.
\]
This means that system $\mathcal{C}_{2(b)}$ is contextual if and
only if
\begin{multline*}
\delta\left(\left(R_{1}^{1},R_{1}^{2}\right),\left(R_{2}^{1},R_{2}^{2}\right)\right)=\min_{\begin{array}{c}
\textnormal{all couplings}\\
\begin{array}{c}
\left(S_{1}^{1},S_{2}^{1},S_{1}^{2},S_{2}^{2}\right)\\
\textnormal{ of system }\mathcal{C}_{2(b)}
\end{array}
\end{array}\textnormal{ }}\left(\Pr\left[S_{1}^{1}\not=S_{1}^{2}\right]+\Pr\left[S_{2}^{1}\not=S_{2}^{2}\right]\right)\\
>\left|a_{1}-a_{2}\right|+\left|b_{1}-b_{2}\right|.
\end{multline*}
Indeed, this inequality indicates that in all couplings either 
\[
\Pr\left[S_{1}^{1}\not=S_{1}^{2}\right]>\delta\left(R_{1}^{1},R_{1}^{2}\right),
\]
 or 
\[
\Pr\left[S_{2}^{1}\not=S_{2}^{2}\right]>\delta\left(R_{2}^{1},R_{2}^{2}\right),
\]
 or both. The intuition remains the same as above: the contexts ``force''
the same-content variables to be more dissimilar than they are in
isolation. The difference 
\[
\delta\left(\left(R_{1}^{1},R_{1}^{2}\right),\left(R_{2}^{1},R_{2}^{2}\right)\right)-\delta\left(R_{1}^{1},R_{1}^{2}\right)-\delta\left(R_{2}^{1},R_{2}^{2}\right)
\]
is a natural (although by far not the only) measure of the degree
of contextuality.\footnote{For other measures of contextuality, see Refs. \cite{CD2020Impossible,DKC2020,DKC2020Erratum,KujDzh2019}} 

\section{Other examples}

The system $\mathcal{C}_{2(b)}$ of the previous section, with the
within-context distributions (\ref{eq:context1 of system1-1}) and
(\ref{eq:context2 of system1-1}), is not a toy example, despite its
simplicity. Except for the specific choice of the questions, it describes
an empirical situation one sees in polls of public opinion, with two
questions asked in one order of a large group of participants, and
the same two questions asked in the other order of another large group
of participants \cite{Moore,WangBusemeyer2013}. 

In quantum physics, system of the $\mathcal{C}_{2}$-type can describe
the outcomes of successive measurements of two spins along two directions,
encoded $1$ and $2$, in the same spin-$\nicefrac{1}{2}$ particle
(e.g., electron). Without getting into details, in such an experiment
the spin-$\nicefrac{1}{2}$ particles are \emph{prepared} in one and
the same \emph{quantum state}, and then subjected to two measurements
in one of the two orders. Each measurement results in one of two outcomes,
spin up ($+1$) or spin down ($-1$). 
\begin{equation}
\begin{array}{|c|c||c|}
\hline R_{1}^{1} & R_{2}^{1} & c^{1}=q_{1}\rightarrow q_{2}\\
\hline R_{1}^{2} & R_{2}^{2} & c^{2}=q_{2}\rightarrow q_{1}\\
\hline\hline q_{1}=\textnormal{"is spin in direction 1 up?"} & q_{2}=\textnormal{"is spin in direction 2 up?"} & \textnormal{system }\mathcal{C}_{2(c)}
\\\hline \end{array}\:.
\end{equation}
The computations in accordance with the standard quantum-mechanical
rules yield the following two results \cite{DzhafarovZhangKujala(2015)Isthere}.
First, the system is inconsistently connected, i.e. generally the
probability of spin-up in a given direction depends on whether it
is measured first or second,
\[
\Pr\left[R_{1}^{1}=+1\right]\not=\Pr\left[R_{1}^{2}=+1\right]\textnormal{ and }\Pr\left[R_{2}^{1}=+1\right]\not=\Pr\left[R_{2}^{2}=+1\right].
\]
Second, the system is noncontextual,\footnote{For those familiar with CbD, this follows from the fact the expected
values $\left\langle R_{1}^{1}R_{2}^{1}\right\rangle $ and $\left\langle R_{1}^{2}R_{2}^{2}\right\rangle $
are always equal to each other, whereas the criterion for contextuality
of a cyclic system \cite{KujalaDzhafarov(2016)Proof}, when specialized
to $n=2$, is $\left|\left\langle R_{1}^{1}R_{2}^{1}\right\rangle -\left\langle R_{1}^{2}R_{2}^{2}\right\rangle \right|>\left|\left\langle R_{1}^{1}\right\rangle -\left\langle R_{1}^{2}\right\rangle \right|+\left|\left\langle R_{2}^{1}\right\rangle -\left\langle R_{2}^{2}\right\rangle \right|.$} i.e., it is always the case that 
\[
\delta\left(\left(R_{1}^{1},R_{1}^{2}\right),\left(R_{2}^{1},R_{2}^{2}\right)\right)\leq\delta\left(R_{1}^{1},R_{1}^{2}\right)+\delta\left(R_{2}^{1},R_{2}^{2}\right).
\]
As we see, systems of the $\mathcal{C}_{2}$-type may be of interest
in both physics and behavioral studies. 

\begin{figure}
\begin{centering}
\includegraphics[scale=0.5]{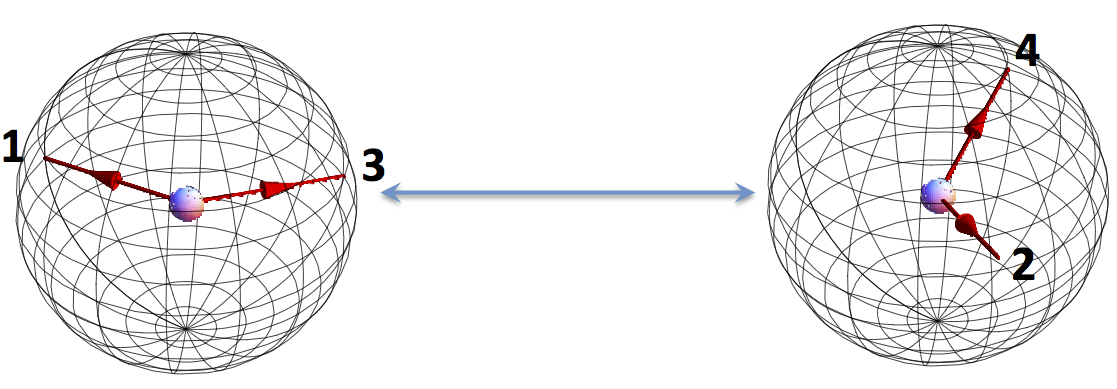}
\par\end{centering}
\caption{\label{fig:EPR/B}A schematic representation of the EPR/Bohm experimental
set up. Explanations in the text.}

\end{figure}

However, in both these fields, the origins of the research of what
we now call contextuality are dated back to another cyclic system,
in which the arrangement shown in Figure \ref{fig:A-cyclic-system}
specializes to
\[
\vcenter{\xymatrix@C=1cm{ & {\scriptstyle R_{1}^{1}}\ar@{-}[rr]_{\textnormal{context }c^{1}} &  & {\scriptstyle R_{2}^{1}}\ar@{.}[dr]_{\textnormal{content }q_{2}}\\
{\scriptstyle R_{1}^{4}}\ar@{.}[ur]_{\textnormal{content }q_{1}} &  &  &  & {\scriptstyle R_{2}^{2}}\ar@{-}[dd]_{\textnormal{context }c^{2}}\\
\\
{\scriptstyle R_{4}^{4}}\ar@{-}[uu]_{\textnormal{context }c^{4}} &  &  &  & {\scriptstyle R_{3}^{2}}\ar@{.}[dl]_{\textnormal{content }q_{3}}\\
 & {\scriptstyle R_{4}^{3}}\ar@{.}[lu]_{\textnormal{content }q_{4}} &  & {\scriptstyle R_{3}^{3}}\ar@{-}[ll]_{\textnormal{context }c^{3}}
}
}.
\]
Figure \ref{fig:EPR/B} illustrates the empirical situation described
by this system, and the first for which contextuality was mathematically
established \cite{Bell1966,CH1974,CHSH1969,BohmAharonov1957,Fine1982JMP}.
Two spin-$\nicefrac{1}{2}$ particles are prepared in a special quantum
state making them \emph{entangled}, and they move away from each other.
The ``left'' particle's spin is measured along one of the two directions
(encoded $1$ and $3$) by someone we will call Zora, and simultaneously
the ``right'' particle's spin is measured along one of the two directions
(encoded $2$ and $4$) by a Nico.\footnote{For no deep reason, I decided to deviate from the established tradition
to call the imaginary performers of the measurements in this task
Alice and Bob.} The outcomes of the measurements are spin-up or spin-down, and each
random variable $R_{i}^{j}$ answers the question 
\[
q_{i}:\textnormal{is the spin in direction }i\textnormal{ up? }(i=1,2,3,4).
\]
 
In the form of a content-context matrix the system can be presented
as
\begin{equation}
\begin{array}{|c|c|c|c||c|}
\hline R_{1}^{1} & R_{2}^{1} &  &  & c^{1}\\
\hline  & R_{2}^{2} & R_{3}^{2} &  & c^{2}\\
\hline  &  & R_{3}^{3} & R_{4}^{3} & c^{3}\\
\hline R_{1}^{4} &  &  & R_{4}^{4} & c^{4}\\
\hline\hline \begin{array}{c}
q_{1}\\
\textnormal{(Zora's 1)}
\end{array} & \begin{array}{c}
q_{2}\\
\textnormal{(Nico's 2)}
\end{array} & \begin{array}{c}
q_{3}\\
\textnormal{(Zora's 3)}
\end{array} & \begin{array}{c}
q_{4}\\
\textnormal{(Nico's 4)}
\end{array} & \mathcal{C}_{4(a)}
\\\hline \end{array}\:.\label{eq:C4}
\end{equation}
The measurements by Zora and Nico are made simultaneously, or at least
close enough in time so that no signal about Zora's choice of a direction
can reach Nico before he makes his measurement, and vice versa. Because
of this, the system is consistently connected,
\[
R_{i}^{j}\overset{dist}{=}R_{i}^{j'}
\]
for any content $q_{i}$ and two contexts $c^{j}$ and $c^{j'}$ in
which $q_{i}$ is measured. Following the logic of contextuality analysis,
we first establish that (because of the consistent connectedness)
\[
\delta\left(R_{1}^{1},R_{1}^{4}\right)=\delta\left(R_{2}^{1},R_{2}^{2}\right)=\delta\left(R_{3}^{2},R_{3}^{3}\right)=\delta\left(R_{4}^{3},R_{4}^{4}\right)=0.
\]
Then we compute 
\begin{multline*}
\delta\left(\left(R_{1}^{1},R_{1}^{4}\right),\left(R_{2}^{1},R_{2}^{2}\right),\left(R_{3}^{2},R_{3}^{3}\right),\left(R_{4}^{3},R_{4}^{4}\right)\right)\\
\\
=\min_{\begin{array}{c}
\textnormal{all couplings}\\
\begin{array}{c}
\left(S_{1}^{1},S_{4}^{1},S_{1}^{2},S_{2}^{2},S_{2}^{3},S_{3}^{3},S_{3}^{4},S_{4}^{4}\right)\\
\textnormal{ of system }\mathcal{C}_{4(a)}
\end{array}
\end{array}\textnormal{ }}\left(\Pr\left[S_{1}^{1}\not=S_{1}^{4}\right]+\Pr\left[S_{2}^{1}\not=S_{2}^{2}\right]+\Pr\left[S_{3}^{2}\not=S_{3}^{3}\right]+\Pr\left[S_{4}^{3}\not=S_{4}^{4}\right]\right).
\end{multline*}
The system is noncontextual if and only if this quantity is zero.
As it turns out (and this is what was established by John Bell in
his celebrated papers in the 1960s, \cite{Bell1964,Bell1966}), the
directions $1,2,3,4$ can be chosen so that, by the laws of quantum
mechanics, this quantity is greater than zero, making the system contextual.

In psychology, systems of the same $\mathcal{C}_{4}$-type have been
of interest as representing the following empirical situation \cite{DzhafarovKujala(2014)TopicsCogSci,Dzhafarov(2003),DzhafarovKujala(2010),DzhafarovKujala(2016)NHMP,Sternberg1969,Townsend1984,KujalaDzhafarov(2008),ZhangDzhafarov(2015)}.
Consider two variables having two values each, that can be manipulated
in an experiment. Think, e.g., of a briefly presented visual object
that can have one of two colors (red or green) and one of two shapes
(square or oval), combined in the $2\times2$ ways. In the experiment,
an observer responds to the object by answering two Yes-No questions:
``is the object red?'' and ``is the object square?''. If we simply
identify these questions with contents, the resulting system of random
variables looks like this:
\begin{equation}
\begin{array}{|c|c||l|}
\hline R_{1}^{1} & R_{2}^{1} & c^{1}:\textnormal{red and oval}\\
\hline R_{1}^{2} & R_{2}^{2} & c^{2}:\textnormal{green and oval}\\
\hline R_{1}^{3} & R_{2}^{3} & c^{3}:\textnormal{red and square}\\
\hline R_{1}^{4} & R_{2}^{4} & c^{4}:\textnormal{green and square}\\
\hline\hline q_{1}:\textnormal{red?} & q_{2}:\textnormal{square?} & \mathcal{R}
\\\hline \end{array}\:,\label{eq:system shape color 1}
\end{equation}
with the contexts describing the object being presented, and the contents
the questions asked. 

Although possible, this is not, however, an especially interesting
way of conceptualizing the situation. It is more informative to describe
the contents of the random variables as color and shape responses
to the color and shape of the visual stimuli, respectively: 
\[
\begin{array}{|l|}
\hline q_{1}:\textnormal{does this red object appear red?}\\
\hline q_{2}:\textnormal{does this square object appear square?}\\
\hline q_{3}:\textnormal{does this green object appear red?}\\
\hline q_{4}:\textnormal{does this oval object appear square?}
\\\hline \end{array}
\]
With the contexts remaining as they are in system (\ref{eq:system shape color 1}),
the experiment is now represented by a system of the $\mathcal{C}_{4}$-type:
\[
\begin{array}{|c|c|c|c||c|}
\hline R_{1}^{1} & R_{2}^{1} &  &  & c^{1}\\
\hline  & R_{2}^{2} & R_{3}^{2} &  & c^{2}\\
\hline  &  & R_{3}^{3} & R_{4}^{3} & c^{3}\\
\hline R_{1}^{4} &  &  & R_{4}^{4} & c^{4}\\
\hline\hline \begin{array}{c}
q_{1}\\
\textnormal{(red)}
\end{array} & \begin{array}{c}
q_{2}\\
\textnormal{(square)}
\end{array} & \begin{array}{c}
q_{3}\\
\textnormal{(green)}
\end{array} & \begin{array}{c}
q_{4}\\
\textnormal{(oval)}
\end{array} & \mathcal{C}_{4(b)}
\\\hline \end{array}
\]
Compared to system $\mathcal{C}_{4(c)}$ in (\ref{eq:C4}), the physical
situation described by $\mathcal{C}_{4(b)}$ is, of course, very different:
e.g., instead of $R_{1}^{j}$ and $R_{3}^{j}$ being outcomes of spin
measurements by Zora along two different directions, these random
variables represent now responses to the color question when the color
is red and when it is green, respectively. However, the logic of the
contextuality analysis does not change. If this system turns out to
be consistently connected and noncontextual, the interpretation of
this in psychology is that the judgment of color is \emph{selectively
influenced} \emph{by object's color} (irrespective of its shape),
and the judgment of shape is \emph{selectively influenced by object's
shape} (irrespective of its color). Deviations from this pattern of
selective influences, whether in the form of inconsistent connectedness
or contextuality, or both,\footnote{System $\mathcal{C}_{4(d)}$ is almost certainly inconsistently connected
(guessing of an imaginary experiment based on the results of many
real ones). } provide an interesting way of classifying (and quantifying) the ways
object's color may influence one's judgment of its shape and vice
versa. 

\section{What if the system is deterministic?}

A \emph{deterministic quantity} $r$ is a special case of a random
variable: it is a random variable $R$ that attains the value $r$
with probability 1:
\[
\Pr\left[R=r\right]=1.
\]
It is convenient to present this as 
\[
R\equiv r.
\]
A \emph{deterministic system} is one containing only deterministic
variables. For instance,
\begin{equation}
\begin{array}{|c|c|c|c|c||cc|}
\hline r_{1}^{1} & r_{2}^{1} &  & r_{4}^{1} &  &  & c^{1}\\
\hline r_{1}^{2} &  & r_{3}^{2} &  &  &  & c^{2}\\
\hline  & r_{2}^{3} & r_{3}^{3} & r_{4}^{3} & r_{5}^{3} &  & c^{3}\\
\hline  &  & r_{3}^{4} &  & r_{5}^{4} &  & c^{4}\\
\hline\hline  &  &  &  &  &  & \\
q_{1} & q_{2} & q_{3} & q_{4} & q_{5} &  & \mathcal{D}
\\\hline \end{array}\label{eq:deterministic}
\end{equation}
is a deterministic systems in which $r_{i}^{j}$ represents a random
variable $R_{i}^{j}\equiv r_{i}^{j}$. The system can be consistently
connected (if the value of $r_{i}^{j}$ does not depend on $j$) or
inconsistently connected (otherwise). 

It is easy to see, however, that a deterministic system is always
noncontextual.\footnote{This fact was first mentioned to me years ago by Matt Jones of the
University of Colorado.} Indeed, any two content-sharing $R_{i}^{j}\equiv r_{i}^{j}$ and
$R_{i}^{j'}\equiv r_{i}^{j'}$ in this system have a single coupling
($S_{i}^{j}\equiv r_{i}^{j}$,$S_{i}^{j}\equiv r_{i}^{j}$), consisting
of the same deterministic quantities but considered jointly distributed.\footnote{There is a subtlety here, first pointed out to me by Janne Kujala
of Turku University. If $R_{i}^{j}\equiv r_{i}^{j}$ and $R_{i}^{j'}\equiv r_{i}^{j'}$,
one may be tempted to say that the joint event $\left(R_{i}^{j}\equiv r_{i}^{j},R_{i}^{j'}\equiv r_{i}^{j'}\right)$
has the probability one, and this would create an exception from the
principle that random variables in different contexts are not jointly
distributed. This is wrong, however, because $\left(R_{i}^{j}\equiv r_{i}^{j},R_{i}^{j'}\equiv r_{i}^{j'}\right)$
can only be thought of counterfactually, as it involves mutually exclusive
contexts. In fact, the only justification (or, better put, excuse)
for the intuition that $\left(R_{i}^{j}\equiv r_{i}^{j},R_{i}^{j'}\equiv r_{i}^{j'}\right)$
is a meaningful joint event is that $R_{i}^{j}\equiv r_{i}^{j}$ and
$R_{i}^{j'}\equiv r_{i}^{j'}$ have a single coupling, and in this
coupling $\Pr\left[S_{i}^{j}\equiv r_{i}^{j},S_{i}^{j}\equiv r_{i}^{j}\right]=1$.
More generally, use of couplings is a rigorous way of dealing with
counterfactuals \cite{Dzhafarov2019}.} It follows that
\[
\delta\left(r_{i}^{j},r_{i}^{j'}\right)=\left\{ \begin{array}{ccc}
1 & \textnormal{if} & r_{i}^{j}\not=r_{i}^{j'}\\
0 & \textnormal{if} & r_{i}^{j}=r_{i}^{j'}
\end{array}\right..
\]
The entire deterministic system in (\ref{eq:deterministic}) also
has a single coupling, one containing the same deterministic quantities
as the system itself, but considered jointly distributed. Clearly,
the subcoupling $\left(S_{i}^{j}\equiv r_{i}^{j},S_{i}^{j'}\equiv r_{i}^{j'}\right)$
extracted from this coupling is precisely the same as the coupling
of $R_{i}^{j}\equiv r_{i}^{j}$ and $R_{i}^{j'}\equiv r_{i}^{j'}$
taken in isolation, and
\[
\delta\left(\left\{ \left(r_{i}^{j},r_{i}^{j'}\right):\textnormal{all such pairs}\right\} \right)=\sum_{\textnormal{all such pairs}}\delta\left(r_{i}^{j},r_{i}^{j'}\right).
\]

One might conclude that deterministic systems are of no interest for
contextuality analysis. This is not always true, however. There are
cases when we know that a system is deterministic, but we do not know
which of a set of possible deterministic systems it is, because it
can be any of them. Let us look at this in detail, using as examples
systems consisting of logical truth values of various statements. 

Consider first the following $\mathcal{C}_{4}$-type system:
\begin{equation}
\begin{array}{|c|c|c|c||c|}
\hline R_{1}^{1}\equiv+1 & R_{2}^{1}\equiv-1 &  &  & c^{1}\\
\hline  & R_{2}^{2}\equiv+1 & R_{3}^{2}\equiv-1 &  & c^{2}\\
\hline  &  & R_{3}^{3}\equiv+1 & R_{4}^{3}\equiv-1 & c^{3}\\
\hline R_{1}^{4}\equiv-1 &  &  & R_{4}^{4}\equiv+1 & c^{4}\\
\hline\hline q_{1} & q_{2} & q_{3} & q_{4} & \mathcal{C}_{4(c)}
\\\hline \end{array}\:,\label{eq:ZoraNiconames}
\end{equation}
where $+1$ and $-1$ encode truth values (true and false), and the
contents are the statements
\[
\begin{array}{|l|c|}
\hline q_{1}:\textnormal{"my name is Zora"} & q_{2}:\textnormal{"my name is Nico"}\\
\hline q_{3}:\textnormal{"my name is Max"} & q_{4}:\textnormal{"my name is Alex"}
\\\hline \end{array}\:.
\]
Equivalently, the contents could also be formulated as questions,
``is my name Zora?'' and ``is my name Nico?'', in which case $+1$
and $-1$ would encode answers Yes and No. In the following, however,
I will refer to the $q$'s as statements, and the values of the variables
as truth values. The contexts justifying the truth values in (\ref{eq:ZoraNiconames})
are
\[
\begin{array}{|l|c|}
\hline c^{1}:\textnormal{the statements are made by Zora} & c^{2}:\textnormal{the statements are made by Nico}\\
\hline c^{3}:\textnormal{the statements are made by Max} & c^{4}:\textnormal{the statements are made by Alex}
\\\hline \end{array}\:.
\]
This is a situation when the truth values are determined uniquely,
the system is deterministic, and consequently it is noncontextual
(even though context-dependence in it is salient in the form of inconsistent
connectedness).

Consider next another system of the $C_{4}$-type, 
\[
\begin{array}{|c|c|c|c||c|}
\hline R_{1}^{1} & R_{2}^{1} &  &  & c^{1}\\
\hline  & R_{2}^{2} & R_{3}^{2} &  & c^{2}\\
\hline  &  & R_{3}^{3} & R_{4}^{3} & c^{3}\\
\hline R_{1}^{4} &  &  & R_{4}^{4} & c^{4}\\
\hline\hline q_{1}:\textnormal{"}q_{2}\textnormal{ is true"} & q_{2}:\textnormal{"}q_{3}\textnormal{ is true"} & q_{3}:\textnormal{"}q_{4}\textnormal{ is true"} & q_{4}:\textnormal{"}q_{1}\textnormal{ is false"} & \mathcal{C}_{4(d)}
\\\hline \end{array}\:,
\]
with contents/statements of a very different kind, and the contexts
which here (at least provisionally) can simply be defined by which
statements they include: $c^{1}$ includes $\left(q_{1},q_{2}\right)$,
$c^{2}$ includes $\left(q_{2},q_{3}\right)$, etc. 

One can recognize here a formalization of the quadripartite version
of the Liar antinomy: one can begin with any statement, say $q_{3}$,
assume it is true, conclude that then $q_{4}$ is true, then $q_{1}$
is false, then $q_{2}$ is false, and then $q_{3}$ is false; and
if one assumes that $q_{3}$ is false, then by the analogous chain
of assignments one arrives to $q_{3}$ being true. There is no consistent
assignment of truth values in this system. In the language of CbD,
the truth values of the statements in $\mathcal{C}_{4(d)}$ can only
be described by an inconsistently connected deterministic system. 

We come to the main issue now: $\mathcal{C}_{4(d)}$ is certainly
a deterministic system (because truth values of statements within
a context are fixed), but which deterministic system is it? There
are 16 possible ways of filling this system with truth values:

\[
\begin{array}{c}
\begin{array}{|c|c|c|c|}
\hline +1 & +1 &  & \\
\hline  & +1 & +1 & \\
\hline  &  & +1 & +1\\
\hline -1 &  &  & +1
\\\hline \end{array}\:,\quad\begin{array}{|c|c|c|c|}
\hline +1 & +1 &  & \\
\hline  & +1 & +1 & \\
\hline  &  & -1 & -1\\
\hline +1 &  &  & -1
\\\hline \end{array}\:,\\
\\
\begin{array}{|c|c|c|c|}
\hline +1 & +1 &  & \\
\hline  & -1 & -1 & \\
\hline  &  & +1 & +1\\
\hline -1 &  &  & +1
\\\hline \end{array}\:,\quad\begin{array}{|c|c|c|c|}
\hline +1 & +1 &  & \\
\hline  & +1 & +1 & \\
\hline  &  & -1 & -1\:\\
\hline +1 &  &  & -1
\\\hline \end{array},\\
\textnormal{etc.}
\end{array}
\]
The only constraint in generating these systems is that 
\begin{enumerate}
\item in the first three contexts (rows) the truth values of the two variables
coincide (because the first statement in them says that the second
one is true, and the second one does not refer to the first one);
\item in context $c^{4}$ (the last row) the truth values of the two variables
are opposite (because $q_{4}$ says that $q_{1}$ is false, and $q_{1}$
does not refer to $q_{4}$).
\end{enumerate}
We see that although random variability in $\mathcal{C}_{4(d)}$ is
absent, we have in its place \emph{epistemic uncertainty}. This opens
the possibility of attaching epistemic (Bayesian) probabilities to
the 16 possible deterministic variants of $\mathcal{C}_{4(d)}$, and
obtaining as a result a \emph{system of epistemic random variables}.
Mathematically, such a variable is treated in precisely the same way
as an ordinary (``frequentist'') random variable. For instance,
we can say that an epistemic variable $R$ can have values $+1$ and
$-1$ with Bayesian probabilities $p$ and $1-p$. This means that
$R$ in fact is a deterministic quantity that can be either $+1$
or $-1$, and the degree of rational belief that $R$ is $+1$ (given
what we know of it) is $p$. In all computational respects, however,
$R$ is treated as if it was a variable that sometimes can be $+1$
and sometimes $-1$.

If we choose equal weights for all 16 deterministic variants of $\mathcal{C}_{4(d)}$
(simply because we have no rational grounds for preferring some of
them to others), the resulting system will have the following Bayesian
distributions:
\begin{equation}
\begin{array}{c|c|c|c}
\begin{array}{c}
\textnormal{context }c^{i},\\
i=1,2,3
\end{array} & R_{i+1}^{i}=+1 & R_{i+1}^{i}=-1\\
\hline R_{i}^{i}=+1 & \nicefrac{1}{2} & 0 & \nicefrac{1}{2}\\
\hline R_{i}^{i}=-1 & 0 & \nicefrac{1}{2} & \nicefrac{1}{2}\\
\hline  & \nicefrac{1}{2} & \nicefrac{1}{2}
\end{array}\label{eq:Liar1}
\end{equation}
and

\begin{equation}
\begin{array}{c|c|c|c}
\textnormal{context }c^{4} & R_{1}^{4}=+1 & R_{1}^{4}=-1\\
\hline R_{4}^{4}=+1 & 0 & \nicefrac{1}{2} & \nicefrac{1}{2}\\
\hline R_{4}^{4}=-1 & \nicefrac{1}{2} & 0 & \nicefrac{1}{2}\\
\hline  & \nicefrac{1}{2} & \nicefrac{1}{2}
\end{array}\:.\label{eq:Liar2}
\end{equation}
This system is clearly contextual. Indeed, since it is consistently
connected,
\begin{equation}
\delta\left(R_{1}^{1},R_{1}^{4}\right)=\delta\left(R_{2}^{1},R_{2}^{2}\right)=\delta\left(R_{3}^{2},R_{3}^{3}\right)=\delta\left(R_{4}^{3},R_{4}^{4}\right)=0.\label{eq:deltasLiar}
\end{equation}
At the same time,
\begin{multline}
\delta\left(\left(R_{1}^{1},R_{1}^{4}\right),\left(R_{2}^{1},R_{2}^{2}\right),\left(R_{3}^{2},R_{3}^{3}\right),\left(R_{4}^{3},R_{4}^{4}\right)\right)\\
\\
=\min_{\begin{array}{c}
\textnormal{all couplings}\\
\begin{array}{c}
\left(S_{1}^{1},S_{4}^{1},S_{1}^{2},S_{2}^{2},S_{2}^{3},S_{3}^{3},S_{3}^{4},S_{4}^{4}\right)\\
\textnormal{ of system }\mathcal{C}_{4(d)}
\end{array}
\end{array}\textnormal{ }}\left(\Pr\left[S_{1}^{1}\not=S_{1}^{4}\right]+\Pr\left[S_{2}^{1}\not=S_{2}^{2}\right]+\Pr\left[S_{3}^{2}\not=S_{3}^{3}\right]+\Pr\left[S_{4}^{3}\not=S_{4}^{4}\right]\right)>0.\label{eq:deltaLiar}
\end{multline}
This is easy to see. This quantity could be zero only if, in some
coupling of $\mathcal{C}_{4(a)}$, the equalities in the first row
below all held with probability 1:
\[
\xymatrix@C=1cm{S_{1}^{4}=S_{1}^{1}\ar[d] & S_{2}^{1}=S_{2}^{2}\ar[d] & S_{3}^{2}=S_{3}^{3}\ar[d] & S_{4}^{3}=S_{4}^{4}\ar[d]\\
S_{1}^{1}=S_{2}^{1}\ar[ru] & S_{2}^{2}=S_{3}^{2}\ar[ru] & S_{3}^{3}=S_{4}^{3}\ar[ru] & S_{4}^{4}\not=S_{1}^{4}
}
\:.
\]
But in any coupling of $\mathcal{C}_{4(a)}$, the equalities in the
second row also hold with probability 1, because they copy (\ref{eq:Liar1})
and (\ref{eq:Liar2}). Reading now all the equalities above from left
to right along the arrows as a chain 
\[
S_{1}^{4}=S_{1}^{1}=S_{2}^{1}=S_{2}^{2}=\ldots,
\]
one arrives at a contradiction
\[
S_{1}^{4}\not=S_{1}^{4}.
\]
In essence, this is the same reasoning as that establishing the unremovable
contraction in the Liar antinomy. However, this time it merely serves
the purpose of establishing that our system is contextual. In fact,
the degree of contextuality here, computed as the difference between
(\ref{eq:deltaLiar}) and the (zero) sum of the deltas in (\ref{eq:deltasLiar}),
is maximal among all possible systems of the $\mathcal{C}_{4}$-type.

\begin{figure}
\begin{centering}
\includegraphics[viewport=0bp 200bp 800bp 1035bp,scale=0.25]{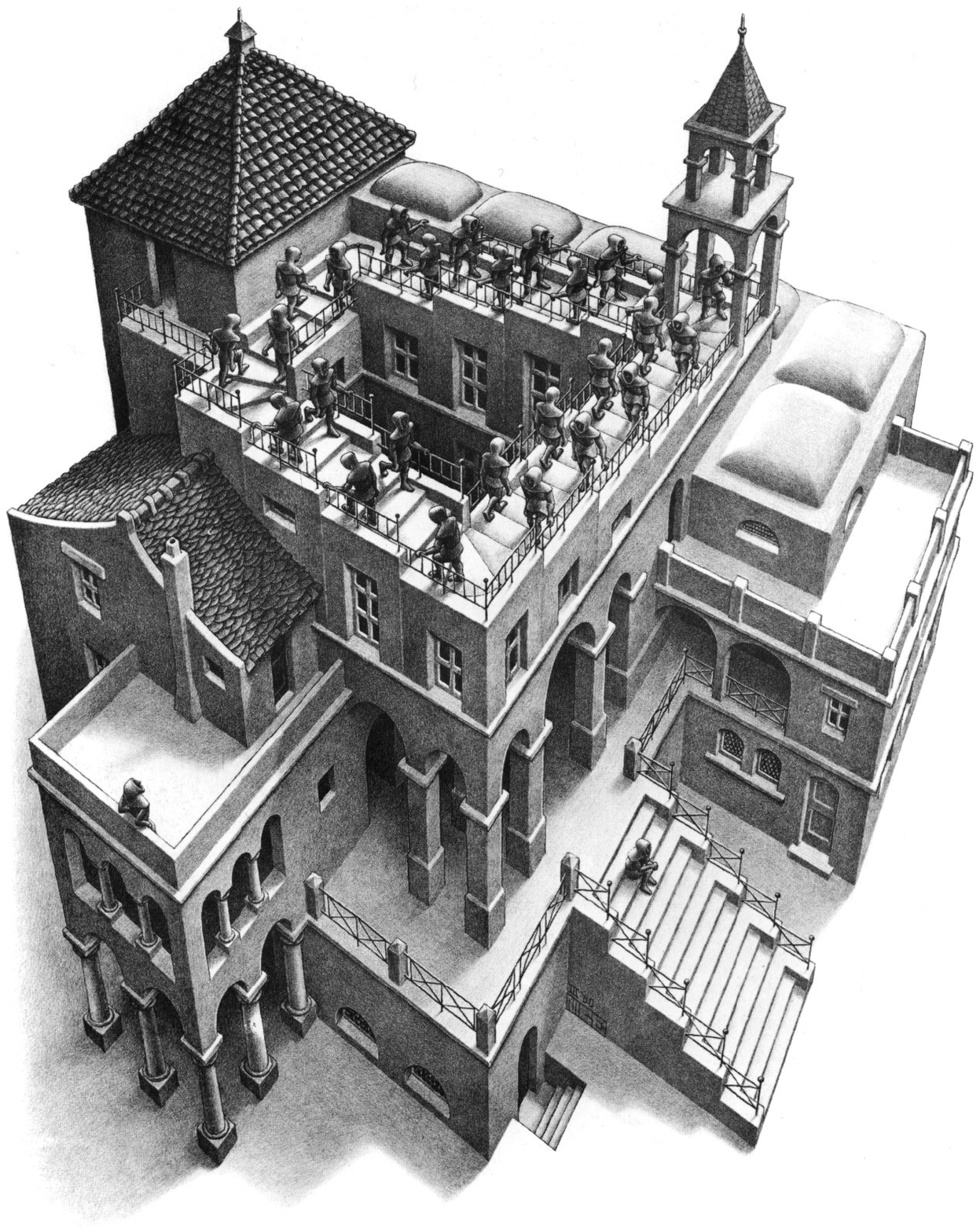}$\begin{array}{c}
\begin{array}{|c|c|c|c||c|}
\hline R_{1}^{1} & R_{2}^{1} &  &  & c^{1}\\
\hline  & R_{2}^{2} & R_{3}^{2} &  & c^{2}\\
\hline  &  & R_{3}^{3} & R_{4}^{3} & c^{3}\\
\hline R_{1}^{4} &  &  & R_{4}^{4} & c^{4}\\
\hline R_{1}^{5} & R_{1}^{5} & R_{1}^{5} & R_{1}^{5} & c^{5}\\
\hline\hline q_{1} & q_{2} & q_{3} & q_{4} & \mathcal{E}
\\\hline \end{array}\end{array}$
\par\end{centering}
\caption{\label{fig:Ascending-Descending}``Ascending-Descending'' by M.
C. Escher. The four flights of stairs are enumerated $q_{1},q_{2},q_{3},q_{4}$.
The epistemic random variables have values \emph{ascending} and \emph{descending},\emph{
}and in each of the the first four contexts they are perfectly correlated.
The fifth context is a mixture of the quadruples of values precisely
two of which are ascending (so that travelers always end up in the
same place from where they started). The resulting epistemic system
is contextual \cite{CD2020Impossible}.}
\end{figure}

We could use other multipartite versions of the Liar paradox, with
three or five or any number of statements, all leading to the same
outcome. A special mention is needed of the bipartite version. In
this system it is no longer possible to define the contexts simply
by the contents of the variables they include. Instead we once again
need to use the order of the contents, this time interpreted as the
direction of inference: $q\rightarrow q'$ means that we assign truth
values to $q$ and infer the corresponding truth values for $q'$.\footnote{The interpretation of contexts in terms of the direction of inference
is the right one also in systems with larger number of statements.
It is merely a coincidence that for $n>2$ in the systems depicting
the $n$-partite Liar paradox the direction of inference in a context
is uniquely determined by the pairs of contents involved in this context.} The resulting system is

\[
\begin{array}{|c|c||c|}
\hline R_{1}^{1} & R_{2}^{1} & c^{1}:q_{1}\rightarrow q_{2}\\
\hline R_{1}^{2} & R_{2}^{2} & c^{2}:q_{2}\rightarrow q_{1}\\
\hline\hline q_{1}=\textnormal{"}q_{2}\textnormal{ is true"} & q_{2}:\textnormal{"}q_{1}\textnormal{ is false"} & \mathcal{C}_{2(d)}
\\\hline \end{array}\:,
\]
with four possible deterministic variants:
\[
\begin{array}{c}
\begin{array}{|c|c|}
\hline +1 & +1\\
\hline -1 & +1
\\\hline \end{array}\:,\quad\begin{array}{|c|c|}
\hline +1 & +1\\
\hline +1 & -1
\\\hline \end{array}\:,\\
\\
\begin{array}{|c|c|}
\hline -1 & -1\\
\hline -1 & +1
\\\hline \end{array}\:,\quad\begin{array}{|c|c|}
\hline +1 & +1\\
\hline +1 & -1
\\\hline \end{array}\:.
\end{array}
\]
Mixing them with equal epistemic probabilities creates a consistently
connected and highly contextual system (maximally contextual among
all cyclic systems of rank 2).

Logical paradoxes are not, of course, the only application of contextuality
analysis with epistemic random variables. It seems that many ``strange''
or ``paradoxical'' situations can be converted into contextual epistemic
systems \cite{CD2020Impossible,Dzhafarov2020}. Among other applications
are such objects as the Penroses' ``impossible figures'' and M.
C. Escher's pictures (as in Figure \ref{fig:Ascending-Descending}). 

\section{The right to ignore (or not to)}

I will mention now some aspects of the Contextuality-by-Default theory
(CbD) that seem to pose difficulties for understanding. Questions
about them are being asked often and in spite of having been repeatedly
addressed in published literature. 

The most basic aspect of CbD is double indexation of the random variables.
The response to a given question $q$ is a random variable $R_{q}^{c}$
whose identity is determined not only by $q$ but also by the context
$c$ in which $q$ is responded to. This looks innocuous enough, but
it puzzles some when a system being analyzed is consistently connected,
i.e. when changing $c$ in $R_{q}^{c}$ does not change the distribution.
And the puzzlement may increase when our knowledge tells us there
is no possible way in which different contexts $c$ can differently
influence the random variables $R_{q}^{c}$. 

Consider again the system $\mathcal{C}_{4(a)}$ in (\ref{eq:C4}),
from which we date contextuality studies. I reproduce it here for
the reader's convenience:

\[
\begin{array}{|c|c|c|c||c|}
\hline R_{1}^{1} & R_{2}^{1} &  &  & c^{1}\\
\hline  & R_{2}^{2} & R_{3}^{2} &  & c^{2}\\
\hline  &  & R_{3}^{3} & R_{4}^{3} & c^{3}\\
\hline R_{1}^{4} &  &  & R_{4}^{4} & c^{4}\\
\hline\hline \begin{array}{c}
q_{1}\\
\textnormal{(Zora's 1)}
\end{array} & \begin{array}{c}
q_{2}\\
\textnormal{(Nico's 2)}
\end{array} & \begin{array}{c}
q_{3}\\
\textnormal{(Zora's 3)}
\end{array} & \begin{array}{c}
q_{4}\\
\textnormal{(Nico's 4)}
\end{array} & \mathcal{C}_{4(a)}
\\\hline \end{array}\:.
\]
In this system, Nico's choice between directions 2 and $4$ can in
no ways affect Zora's measurements of spin along direction $1$. Nevertheless,
when Nico switches from direction $2$ to $4$, the random variable
describing the outcome of Zora's measurement of spin along direction
$1$ ceases to be $R_{1}^{1}$ and becomes $R_{1}^{4}$. It looks
like Nico has influenced Zora's measurements after all. Isn't it an
example of what Albert Einstein famously called a ``spooky action
at a distance''? 

The answer is, it is not. Nico's choices are undetectable by Zora.
Whether he chooses direction 2 or direction 4, Zora can see no changes
in the statistical properties of what she observes when she measures
spins along direction 1. ``Action'' means information transmitted,
and no information is transmitted from Nico to Zora (and vice versa).
The fact that in at least one of the pairs 
\[
\left\{ R_{1}^{1},R_{1}^{4}\right\} ,\left\{ R_{2}^{1},R_{2}^{2}\right\} ,\left\{ R_{3}^{3},R_{3}^{4}\right\} ,\left\{ R_{4}^{3},R_{4}^{4}\right\} 
\]
the two random variables cannot be viewed as being the same can be
established by neither Zora nor Nico. It can only be established by
a Max who receives the choice of directions and outcomes of measurements
from both Zora and Nico and computes the joint distributions in contexts
$c^{1},c^{2},c^{3},c^{4}$.

An important point here is that compared to Max, Zora does not misunderstand
or miss anything when she sees no difference between $R_{1}^{1}$
and $R_{1}^{4}$ or between $R_{3}^{3}$ and $R_{3}^{4}$. Her understanding
is no less complete or less correct. Zora and Max simply deal with
different systems of random variables. In the same way Max's understanding
is no less complete or less correct than that of an Alex who, in addition
to knowing what Max knows, observes whether solar activity during
the measurements is high or low. In Alex's system, each context of
system $\mathcal{C}_{4(a)}$ is split into two contexts, e.g., $c^{1}$
is replaced with
\[
\begin{array}{|c|c|c|c||c|}
\hline R_{1}^{1,high} & R_{2}^{1,high} &  &  & c^{1,high}\\
\hline R_{1}^{1,low} & R_{2}^{1,low} &  &  & c^{1,low}\\
\hline\hline q_{1} & q_{2} & q_{3} & q_{4} & \mathcal{C}_{4(a)}/c^{1}\textnormal{ only}
\\\hline \end{array}\:.
\]

In studying a system of random variable one always can ignore any
of the circumstances that do not affect the distributions of the variables.\footnote{This statement can even be extended to ignoring circumstances when
distributions do change (inconsistent connectedness). However, this
issue has more complex ramifications, and we will set it aside.} Or one can choose not to ignore such circumstances, to systematically
record them and make them part of the contexts. If a circumstance
is irrelevant (as it may be in the case of Alex's recording of solar
activity), one will find this out by considering couplings of the
system. Thus, one may establish that the contextuality analysis of
the system does not change if all couplings are constrained by
\[
\Pr\left[S_{i}^{j,high}=S_{i}^{j,low}\right]=1,
\]
for any $R_{i}^{j}$ in the original system $\mathcal{C}_{4(a)}$.
This would mean that $R_{i}^{j,high}$ and $R_{i}^{j,low}$ can be
viewed as being one and the same random variable (assuming, of course,
that solar activity is indeed irrelevant).

This reasoning fully applies to the issue often raised by those who
enjoy shallow paradoxes. If one records values of a random variable
$R$ in, say, chronological order, and simultaneously records the
ordinal positions of these values in the sequence (as part of their
contexts),
\[
\begin{array}{|c|c|c|c|c|}
\hline r_{1} & r_{2} & \ldots & r_{n} & \ldots\\
\hline 1 & 2 & \ldots & n & \ldots
\\\hline \end{array}\:,
\]
would not this transform all these realizations of a single random
variable into pairwise stochastically unrelated random variables 
\[
R^{1},R^{2},\ldots,R^{n},\ldots
\]
with a single realization each? The answer is yes, if one so wishes
(one may also choose to ignore the ordinal positions of the observations
altogether), but then a standard view is immediately restored when
one considers couplings of these random variables. For instance, the
\emph{iid coupling} (corresponding to the standard statistical concept
of \emph{independent identically distributed variables}) has the structure
\[
\begin{array}{cccccc}
 & R^{1} & R^{2} & \ldots & R^{n} & \ldots\\
S^{1} & \boxed{\overset{}{\:r_{1}^{1}=r_{1}}\:} & r_{2}^{1} & \ldots & r_{n}^{1} & \ldots\\
S^{2} & r_{1}^{2} & \boxed{\overset{}{\:r_{2}^{2}=r_{2}}\:} & \ldots & r_{n}^{2} & \ldots,\\
\vdots & \vdots & \vdots & \ddots & \vdots & \ddots\\
S^{n} & r_{1}^{n} & r_{2}^{n} & \ldots & \boxed{\overset{}{\:r_{n}^{n}=r_{n}}\:} & \ldots\\
\vdots & \vdots & \vdots & \ddots & \vdots & \ddots
\end{array}
\]
where the boxed values are those factually observed, whereas all other
values are independently sampled from the distribution of $R$. More
details are available in Refs. \cite{Dzhafarov(2016),DzhafarovKon2018}.

Finally, does the double-indexation in CbD lend any support to the
holistic view of the universe, the view that ``everything depends
on everything else''? The answer is that the opposite is true, CbD
supports a radically analytic view. First, as we have established,
unless distributions of two given content-sharing variables are found
to be different (which is ubiquitous but not universal) one can ignore
the difference between their contexts, i.e., disregard all other variables
in these contexts. This will redefine the system, but will not be
wrong. Second, the difference in the identity of two content-sharing
variables in different contexts (whether their distributions are the
same or not) involves no change in the colloquial meaning of the word.
The notion of a change implies that something that preserves its identity
(e.g., a moving body) changes some of its properties (e.g., position
in space). However, $R_{2}^{1}$ and $R_{2}^{2}$ (having the same
content in different contexts) are simply different random variables,
stochastically unrelated because they occur in mutually exclusive
contexts. The difference between them is precisely the same as that
between $R_{2}^{1}$ and $R_{1}^{1}$ (different contents in the same
context). By choosing a different question to ask, one switches to
considering another random variable rather than ``changes'' the
previous one. The same happens when one chooses a different context:
one simply switches to considering a different random variable. If
I see Max and then see Alex, it does not mean that Max has changed
into Alex.

The core of these and other problems with understanding CbD, it seems
to me, is in the tendency to view random variables as empirical objects.
They are not. Random variables are our descriptions of empirical objects.
They are part of our knowledge of the world, and the same as any other
knowledge, they can appear, disappear, and be revised as soon as we
adopt a new point of view or gain new evidence.


\begin{thebibliography}{10}
\bibitem{Bell1964}Bell, J. (1964). On the Einstein-Podolsky-Rosen
paradox. Physics 1, 195-200.

\bibitem{Bell1966}Bell, J. (1966). On the problem of hidden variables
in quantum mechanics. Review of Modern Physic 38:447-453.

\bibitem{KochenSpecker(1967)}Kochen, S., Specker, E. P. (1967). The
problem of hidden variables in quantum mechanics. Journal of Mathematics
and Mechanics 17:59--87.

\bibitem{CHSH1969}Clauser, J.F. , Horne, M.A. , Shimony, A., \& Holt,
R.A. (1969). Proposed experiment to test local hidden-variable theories.
Physical Review Letters 23:880-884.

\bibitem{CH1974}Clauser, J.F., Horne, M.A. (1974). Experimental consequences
of objective local theories. Physical Review D 10:526-535.

\bibitem{Fine1982JMP}Fine, A. (1982). Joint distributions, quantum
correlations, and commuting observables. Journal of Mathematical Physics
23:1306-1310.

\bibitem{Cabelloetal(1996)}Cabello, A., Estebaranz, J., \& Garcìa-Alcaine,
G. (1996). Bell-Kochen-Specker theorem: A proof with 18 vectors. Physics
Letters A 212:183-187.

\bibitem{Cabello(2008)}Cabello, A. (2008). Experimentally testable
state-independent quantum contextuality. Physical Review Letters 101,
210401.

\bibitem{KCBS2008}Klyachko, A.A., Can, M.A., Binicioglu, S., \& Shumovsky,
A.S. (2008). A simple test for hidden variables in spin-1 system.
Physical Review Letters 101:020403.

\bibitem{KurzynskiCabelloKaszlikowski(2014)}Kurzynski, P., Cabello,
A., \& Kaszlikowski, D. (2014). Fundamental monogamy relation between
contextuality and nonlocality. Physical Review Letters 112:100401.

\bibitem{Sternberg1969}Sternberg, S. (1969). The discovery of processing
stages: Extensions of Donders' method. In W.G. Koster (Ed.), Attention
and Performance II. Acta Psychologica 30:276--315.

\bibitem{Townsend1984}Townsend, J. T. (1984). Uncovering mental processes
with factorial experiments. Journal of Mathematical Psychology, 28,
363--400.

\bibitem{Dzhafarov(2003)}Dzhafarov, E.N. (2003). Selective influence
through conditional independence. Psychometrika, 68:7-26.

\bibitem{DzhafarovGluhovsky(2006)}Dzhafarov, E.N., \& Gluhovsky,
I. (2006). Notes on selective influence, probabilistic causality,
and probabilistic dimensionality. Journal of Mathematical Psychology,
50:390-401.

\bibitem{KujalaDzhafarov(2008)}Kujala, J.V., Dzhafarov, E.N. (2008).
Testing for selectivity in the dependence of random variables on external
factors. Journal of Mathematical Psychology 52:128--144.

\bibitem{DzhafarovKujala(2010)}Dzhafarov, E.N., Kujala, J.V. (2010).
The joint distribution criterion and the distance tests for selective
probabilistic causality. Frontiers in Psychology: Quantitative Psychology
and Measurement 1:151 doi: 10.3389\slash fpsyg.2010.0015.

\bibitem{DzhafarovKujala(2016)NHMP}Dzhafarov, E.N., \& Kujala, J.V.
(2016). Probability, random variables, and selectivity. In W.Batchelder,
H. Colonius, E.N. Dzhafarov, J. Myung (Eds), The New Handbook of Mathematical
Psychology (pp. 85-150). Cambridge University Press.

\bibitem{ZhangDzhafarov(2015)}Zhang, R. Dzhafarov, E.N. (2015). Noncontextuality
with marginal selectivity in reconstructing mental architectures.
Frontiers in Psychology: Cognition 1:12 doi: 10.3389/fpsyg.2015.00735.

\bibitem{DzhafarovKujala(2012)Selectivity}Dzhafarov, E.N. \& Kujala,
J.V. (2012). Selectivity in probabilistic causality: Where psychology
runs into quantum physics. Journal of Mathematical Psychology 56,
54-63.

\bibitem{DzhafarovKujala(2012)Quantum}Dzhafarov, E.N., Kujala, J.V.
(2012). Quantum entanglement and the issue of selective influences
in psychology: An overview. Lecture Notes in Computer Science 7620:184-195.

\bibitem{DzhafarovKujala(2013)AllPossible}Dzhafarov, E.N., Kujala,
J.V. (2013). All-possible-couplings approach to measuring probabilistic
context. PLoS ONE 8(5):e61712. doi:10.1371/journal.pone.0061712.

\bibitem{DzhafarovKujala(2013)OrderDistance}Dzhafarov, E.N., Kujala,
J.V. (2013). Order-distance and other metric-like functions on jointly
distributed random variables. Proceedings of the American Mathematical
Society 141:3291-3301.

\bibitem{DzhafarovKujala(2014)TopicsCogSci}Dzhafarov, E.N., Kujala,
J.V. (2014). Selective influences, marginal selectivity, and Bell/CHSH
inequalities. Topics in Cognitive Science, 6, 121--128.

\bibitem{DzhafarovKujala(2014)AQualifiesKolmogorov}Dzhafarov, E.N.,
Kujala, J.V. (2014). A qualified Kolmogorovian account of probabilistic
contextuality. Lecture Notes in Computer Science 8369:201-212.

\bibitem{DK2014Scripta-1}Dzhafarov, E.N., Kujala, J.V. (2014). Contextuality
is about identity of random variables. Physica Scripta T163:014009.

\bibitem{DKconversations}Dzhafarov, E.N., Kujala, J.V. (2015). Conversations
on contextuality. In E.N. Dzhafarov, S. Jordan, R. Zhang, V. Cervantes
(Eds). Contextuality from Quantum Physics to Psychology (pp. 1-22).
New Jersey: World Scientific.

\bibitem{KujalaDzhafarov2015}Kujala, J.V., \& Dzhafarov, E.N. (2015).
Probabilistic Contextuality in EPR/Bohm-type systems with signaling
allowed. In E.N. Dzhafarov, S. Jordan, R. Zhang, V. Cervantes (Eds).
Contextuality from Quantum Physics to Psychology (pp. 287-308). New
Jersey: World Scientific.

\bibitem{bacciagaluppi2}Bacciagaluppi, G. (2015). Einsten, Bohm,
and Leggett-Garg. In E.N. Dzhafarov, S. Jordan, R. Zhang, V. Cervantes
(Eds). Contextuality from Quantum Physics to Psychology (pp. 63-76).
New Jersey: World Scientific.

\bibitem{DzhafarovKujalaLarsson(2015)}Dzhafarov, E.N., \& Kujala,
J.V., \& Larsson, J.-Å. (2015). Contextuality in three types of quantum-mechanical
systems. Foundations of Physics 7, 762-782.

\bibitem{DzhafarovZhangKujala(2015)Isthere}Dzhafarov, E.N., Zhang,
R., Kujala, J.V. (2015). Is there contextuality in behavioral and
social systems? Philosophical Transactions of the Royal Society A
374:20150099.

\bibitem{KujalaDzhafarovLar(2015)}Kujala, J.V., Dzhafarov, E.N.,
Larsson, J-Å (2015). Necessary and sufficient conditions for extended
noncontextuality in a broad class of quantum mechanical systems. Physical
Review Letters 115:150401.

\bibitem{DzhafarovKujalaCervantesZhangJones(2016)}Dzhafarov, E.N.,
Kujala, J.V., Cervantes, V.H., Zhang, R., Jones, M. (2016). On contextuality
in behavioral data. Philosophical Transactions of the Royal Society
A 374:20150234.

\bibitem{DzhafarovKujalaCervantes(2016)LNCS}Dzhafarov, E.N., Kujala,
J.V., \& Cervantes, V.H. (2016). Contextuality-by-Default: A brief
overview of ideas, concepts, and terminology. Lecture Notes in Computer
Science 9535:12-23.

\bibitem{DzhafarovKujala(2016)Context-Content}Dzhafarov, E.N., Kujala,
J.V. (2016). Context-content systems of random variables: The contextuality-by-default
theory. Journal of Mathematical Psychology 74:11-33.

\bibitem{deBarrosDzhafarovetal(2016)}de Barros, J.A., Dzhafarov,
E.N., Kujala, J.V., \& Oas, G. (2016). Measuring observable quantum
contextuality. Lecture Notes in Computer Science 9535, 36-47.

\bibitem{KujalaDzhafarov(2016)Proof}Kujala, J.V., Dzhafarov, E.N.
(2016). Proof of a conjecture on contextuality in cyclic systems with
binary variables. Foundations of Physics 46:282-299.

\bibitem{Dzhafarov(2016)}Dzhafarov, E.N. (2016). Stochastic unrelatedness,
couplings, and contextuality. Journal of Mathematical Psychology 75C:34-41.

\bibitem{DzhKuj2017Fortsch-1}Dzhafarov, E.N. and Kujala, J.V. (2017).
Probabilistic foundations of contextuality. Fortschritte der Physik
65:1-11.

\bibitem{DzhafarovKujala(2017)LNCS-2}Dzhafarov, E.N. and Kujala,
J.V. (2017). Contextuality-by-Default 2.0: Systems with binary random
variables. Lecture Notes Computer Sciences 10106:16-32.

\bibitem{CervantesDzhafarov(2016)}Cervantes, V.H., Dzhafarov, E.N.
(2017). Exploration of contextuality in a psychophysical double-detection
experiment. Lecture Notes in Computer Science 10106:182-193.

\bibitem{ZhangDzhafarov(2016)}Zhang, R., Dzhafarov, E.N. (2016).
Testing contextuality in cyclic psychophysical systems of high ranks.
Lecture Notes in Computer Science 10106:151-162.

\bibitem{DzhCerKuj2017}Dzhafarov, E.N., Cervantes, V.H., Kujala,
J.V. (2017). Contextuality in canonical systems of random variables.
Philosophical Transactions of the Royal Society A 375:20160389.

\bibitem{CervantesDzhafarov2017}Cervantes, V.H., Dzhafarov, E.N.
(2017). Advanced analysis of quantum contextuality in a psychophysical
double-detection experiment. Journal of Mathematical Psychology 79:77-84.

\bibitem{Dzh2017Nothing}Dzhafarov, E.N. (2017). Replacing nothing
with something special: Contextuality-by-Default and dummy measurements.
In A. Khrennikov \& T. Bourama (Eds.), Quantum Foundations, Probability
and Information (pp. 39-44). Berlin: Springer.

\bibitem{CervDzh2018}Cervantes, V.H., Dzhafarov, E.N. (2018). Snow
Queen is evil and beautiful: Experimental evidence for probabilistic
contextuality in human choices. Decision 5:193-204.

\bibitem{DzhKuj2018}Dzhafarov, E.N., Kujala, J.V. (2018). Contextuality
analysis of the double slit experiment (with a glimpse into three
slits). Entropy 20:278.

\bibitem{Basievaetall2018}Basieva, I., Cervantes, V.H., Dzhafarov,
E.N., Khrennikov, A. (2019). True contextuality beats direct influences
in human decision making. Journal of Experimental Psychology: General
148, 1925-1937.

\bibitem{CervDzhaf2019}Cervantes, V.H., Dzhafarov, E.N. (2019). True
contextuality in a psychophysical experiment. Journal of Mathematical
Psychology 91, 119-127.

\bibitem{Dzhafarov2019}Dzhafarov, E.N. (2019). On joint distributions,
counterfactual values, and hidden variables in understanding contextuality.
Philosophical Transactions of the Royal Society A 377:20190144.

\bibitem{KujDzh2019}Kujala, J.V., \& Dzhafarov, E.N. (2019). Measures
of contextuality and noncontextuality. Philosophical Transactions
of the Royal Society A 377:20190149. (available as arXiv:1903.07170.)

\bibitem{Jones2019}Jones, M. (2019). Relating causal and probabilistic
approaches to contextuality. Philosophical Transactions of the Royal
Society A, 377, 20190133.

\bibitem{DzhKuj2020}Dzhafarov, E.N., \& Kujala, J.V. (2020). Systems
of random variables and the Free Will Theorem. Physical Review Research
2:043288; doi: 10.1103/PhysRevResearch.2.043288.

\bibitem{DKC2020}Dzhafarov, E.N., Kujala, J.V., \& Cervantes, V.H.
(2020). Contextuality and noncontextuality measures and generalized
Bell inequalities for cyclic systems. Physical Review A 101:042119.

\bibitem{DKC2020Erratum}Dzhafarov, E.N., Kujala, J.V., \& Cervantes,
V.H. (2020). Erratum: Contextuality and noncontextuality measures
and generalized Bell inequalities for cyclic systems {[}Phys. Rev.
A 101, 042119 (2020){]}. Physical Review A 101:069902.

\bibitem{CD2020Impossible}Cervantes, V.H., \& Dzhafarov, E.N. (2020).
Contextuality analysis of impossible figures. Entropy 22, 981; doi:
10.3390/e22090981.

\bibitem{DKC2020DEpistemic}Dzhafarov, E.N., Kujala, J.V., \& Cervantes,
V.H. (2021). Epistemic odds of contextuality in cyclic systems. European
Physics Journal - Special Topics 230:937-940. (available as arXiv:2002.07755.)

\bibitem{Dzhafarov2020}Dzhafarov, E.N. (in press). The Contextuality-by-Default
view of the Sheaf-Theoretic approach to contextuality. To be published
in A. Palmigiano and M. Sadrzadeh (Eds.) Samson Abramsky on Logic
and Structure in Computer Science and Beyond, in series Outstanding
Contributions to Logic. Springer Nature. (available as arXiv:1906.02718.)

\bibitem{Moore}Moore, D.W. (2002). Measuring new types of question-order
effects. Public Opinion Quarterly 66:80-91.

\bibitem{WangBusemeyer2013}Wang, Z., Busemeyer, J.R. (2013). A quantum
question order model supported by empirical tests of an a priori and
precise prediction. Topics in Cognitive Science 5:689--710.

\bibitem{BohmAharonov1957}Bohm, D, Aharonov, Y. (1957). Discussion
of experimental proof for the paradox of Einstein, Rosen and Podolski.
Physical Review 108:1070-1076.

\bibitem{DzhafarovKon2018}Dzhafarov, E.N., Kon, M. (2018). On universality
of classical probability with contextually labeled random variables.
Journal of Mathematical Psychology 85:17-24.
\end{thebibliography}
\end{document}